\documentclass[11pt]{article}
\usepackage[utf8]{inputenc}
\usepackage[english]{babel}
\usepackage{hyperref}
\usepackage[all]{xy}
\usepackage{makeidx}
\usepackage[paper=a4paper,left=30mm,right=20mm,top=25mm,bottom=30mm]{geometry}
\usepackage{tikz}
\usetikzlibrary{fadings,arrows}

\usepackage{amsmath, amssymb, mathrsfs, verbatim, multirow}
\usepackage{pifont}
\usepackage{float}

\newtheorem{thm}{Theorem}
\newtheorem{cor}[thm]{Corollary}
\newtheorem{prop}[thm]{Proposition}
\newtheorem{deft}[thm]{Definition}
\newtheorem{rek}[thm]{Remark}
\newtheorem{conj}{Conjecture}
\let\noi=\noindent

\def\N{\mathbb{N}} 
\def\Z{\mathbb{Z}} 
\def\R{\mathbb{R}} 
\def\C{\mathbb{C}} 
\def\F{\mathbb{F}}

\def\notin{\mbox{$\in$ \hspace{-.8em}/}} 

\makeindex
\title{A note on the Casas-Alvero Conjecture}
\author{Daniel Schaub, Univ Angers, CNRS, LAREMA, SFR MATHSTIC\\F-49000 Angers, France
\\email: daniel.schaub@univ-angers.fr \and
Mark Spivakovsky, Univ Paul Sabatier, CNRS, IMT UMR 5219\\ F-31062 Toulouse, France and\\
CNRS, LaSol UMI 2001, UNAM.
\\email: mark.spivakovsky@math.univ-toulouse.fr }

\begin{document}
\maketitle

\begin{abstract} The Casas--Alvero conjecture predicts that every univariate polynomial $f$ over a field $K$ of characteristic zero having a common factor with each of its derivatives $H_i(f)$ is a power of a linear polynomial. Let $f=x^d+a_1x^{d-1}+\cdots+a_1x \in K[a_1,\ldots,a_{d-1}][x]$ and let $R_i = Res(f,H_i(f))\in K[a_1,\ldots,a_{d-1}]$ be the resultant of $f$ and $H_i(f)$, $i \in \{1,\ldots,d-1\}$. The Casas-Alvero Conjecture is equivalent to saying that $R_1,\ldots,R_{d-1}$ are ``independent'' in a certain sense, namely that the height $ht(R_1,\ldots,R_{d-1})=d-1$ in $K[a_1,\ldots,a_{d-1}]$. In this paper we prove a partial result in this direction: if $i \in \{d-3,d-2,d-1\}$ then
$R_i \notin\sqrt{(R_1,\ldots,\breve{R_i},\ldots,R_{d-1})}$.
\end{abstract}

\section{Introduction}
In the year 2001 Eduardo Casas--Alvero published a paper on higher order polar germs of plane curve singularities \cite{C}. His work on polar germs inspired him to make the following conjecture (according to the testimony of Jos\'e Manuel Aroca, E. Casas communicated the problem orally well before 2001).

Let $K$ be a field, $d$ a strictly positive integer and $f = x^d+a_1x^{d-1}+\cdots+a_{d-1}x+a_d$ a monic univariate polynomial of degree $d$ over $K$. Let
$$
H_i(f) = \binom{d}{i}x^{d-i} + \binom{d-1}{i}a_1x^{d-i-1} + \cdots + \binom{i}{i}a_{d-i}
$$
be the $i$-th Hasse derivative of $f$.

\begin{deft}
The polynomial $f$ is said to be a {\bf Casas--Alvero polynomial} if for each $i\in\{1,\ldots,d-1\}$ it has a non-constant common factor with its $i$-th Hasse derivative $H_i(f)$.
\end{deft}

Note that, by definition, a Casas-Alvero polynomial $f$ has a common root with $H_{d-1}(f)$. In particular, if $\text{char}\ K=0$, it has at least one root $b\in K$, regardless of whether or not $K$ is algebraically closed. Making the change of variables $x\rightsquigarrow x-b$, we may assume that $0$ is a root of $f$, in other words, $a_d=0$. In the sequel, we will always make this assumption without mentioning it explicitly.
\begin{conj} {\bf (Casas--Alvero)}
Assume that $\text{char}\ K=0$. If $f \in K[x]$ is a Casas-Alvero polynomial of degree $d$ with $a_d=0$, then $f(x) =x^d$.
\end{conj}

For $i\in\{1,\ldots,d-1\}$, let $R_i=\text{Res}(f,H_i(f))\in K[a_1,\ldots,a_{d-1}]$ be the resultant of $f$ and $H_i(f)$. The polynomials $f$ and $H_i(f)$ have a common factor if and only if $R_i=0$. Thus $f$ is Casas--Alvero if and only if the point
$(a_1,\ldots,a_{d-1})\in K^{d-1}$ belongs to the algebraic variety $V(R_1,\dots,R_{d-1})\subset K^{d-1}$. In those terms the Conjecture can be reformulated as follows:

\begin{conj} \label{conj:deux}
Let $V=V(R_1,\ldots,R_{d-1}) \subset K^{d-1}$. Then $V =\{0\}$.
\end{conj}

If the field $K$ is algebraically closed then Conjecture \ref{conj:deux} is also eqiuvalent to

\begin{conj}\label{nulstellen} We have
\begin{equation} \label{eq:racine}
\sqrt{(R_1,\ldots,R_{d-1})}=(a_1,\ldots,a_{d-1}) \end{equation} 
or, equivalently,
\begin{equation}
\label{eq:rac2} 
a_i^N \in (R_1,\ldots,R_{d-1}) \text{ for all } i \in \{1,\ldots,d-1\} \text{ and some } N \in \N.
\end{equation}
\end{conj}

For non-algebraically closed fields Conjecture \ref{nulstellen} is a priori stronger than Conjecture \ref{conj:deux}.

\begin{rek}
Let $K \subset K'$ be a field extension. The induced extension
$$
K[a_1,\ldots,a_{d-1}] \subset K'[a_1,\ldots,a_{d-1}]
$$
is faithfully flat. Since the polynomials $R_1,\ldots,R_{d-1}$ have coefficients in $\Z$, (\ref{eq:rac2}) holds in\linebreak
$K[a_1,\ldots,a_{d-1}]$ if and only if it holds in $K'[a_1,\ldots,a_{d-1}]$. Hence the truth of Conjecture \ref{nulstellen} for any given $d$ depends only on the characteristic of $K$ but not on the choice of the field $K$ itself. Because of this, we will take $K=\C$ in the sequel.
\end{rek}

\begin{rek}\label{Grobner}
Formulae (\ref{eq:racine}) and (\ref{eq:rac2}) can be interpreted in terms of Gröbner bases. Namely, (\ref{eq:racine})  and (\ref{eq:rac2}) are equivalent to saying that for any choice of monomial ordering and Gröbner basis $(f_1,\ldots,f_s)$ of
$(R_1,\ldots,R_{d-1})$, after renumbering the $f_j$, the leading monomial of $f_j$ is a power of $a_j$ for
$j\in\{1,\dots,d-1\}$. 
\end{rek}

\begin{rek}\label{rek:gen} Conjecture \ref{nulstellen} and Remark \ref{Grobner} say that, as polynomials in
$K[a_1,\dots,a_{d-1}]$, the resultants $R_1,\dots,R_{d-1}$ are ``independent" in a certain sense.

Each of the following statements is also equivalent to Conjecture \ref{nulstellen}.
\begin{enumerate}
 \item[(a)] For each $i \in \{1,\ldots,d-2\}$, the element $R_{i+1}$ is not a zero divisor modulo $(R_1,\ldots,R_i)$ (in other words, $R_1,\dots,R_{d-1}$ form a regular sequence in $K[a_1,\dots,a_{d-1}]$).
 \item[(b)] For each $i \in \{1,\ldots,d-2\}$, $$R_{i+1} \notin\bigcup\limits_{\mathfrak{p} \in \text{Ass}((R_1,\ldots,R_i))} \mathfrak{p}. $$ where $\text{Ass}((R_1,\ldots,R_i))$ is the set of associated primes of the ideal $(R_1,\ldots,R_i)$.
\end{enumerate}
Moreover, the above statements (a) and (b) are independent of the numbering of the $R_i$; a permutation of the $R_i$ yields equivalent statements.
\end{rek}

\noi{\bf Notation.} We will denote by $(R_1,R_2,\ldots,\breve{R_i},\ldots,R_{d-1})$ the ideal of $K[a_1,\dots,a_{d-1}]$ generated by the set $\{R_1,R_2,\ldots,R_{d-1}\}\setminus\{R_i\}$.
\medskip

The main theorem of this paper is the following partial result in the direction of Conjecture \ref{nulstellen} and statements (a) and (b) of Remark \ref{rek:gen}.

\begin{thm} \label{thm:princ} Take an element $i\in\{d-3,d-2,d-1\}$. We have
$$
R_i \notin\sqrt{(R_1,R_2,\ldots,\breve{R_i},\ldots,R_{d-1})}.
$$
\end{thm}
\medskip

\noindent{\bf Added in press:} In two recent preprints \cite{Gh1} and \cite{Gh2} Soham Ghosh gave a complete proof of the Casas--Alvero conjecture.

\section{Ideals generated by all the resultants but one}

In this section we prove Theorem \ref{thm:princ} after recalling some preliminary results. 

\begin{prop} \label{prop:sturm} Let $f$ be a polynomial of degree $d$ with real roots $\beta_1\le\beta_2\le\ldots\le\beta_d$, counted with multiplicity. Then $H_1(f)$ has real roots $\gamma_1 \le\gamma_2\le\ldots\le\gamma_{d-1}$, counted with multiplicity, where
$\gamma_i\in]\beta_i,\beta_{i+1}[$ if $\beta_i<\beta_{i+1}$ and $\gamma_i=\beta_i$ if $\beta_i=\beta_{i+1}$.
\end{prop}

\noi Proof: Assume that $f$ has $s$ distincts roots $\delta_1<\delta_2<\cdots<\delta_s$ of multiplicities $m_1,\ldots,m_s$, respectively. Then $\delta_j$ is a root of $H_1(f)$ of multiplicity $m_j-1$, where we say that $\delta_j$ is a root of multiplicity 0 if it is not a root of $H_1(f)$. 

By Rolle's theorem, there is at least one root of $H_1(f)$ in each of the $s-1$ open intervals
$]\delta_1,\delta_2[,\ldots,]\delta_{s-1},\delta_s[$. 
\medskip

\noi{\bf Notation.} Let Int$(\beta_i,\beta_{i+1}):=]\beta_i,\beta_{i+1}[$ if $\beta_i<\beta_{i+1}$ and
Int$(\beta_i,\beta_{i+1}):=\{\beta_i\}$ if $\beta_i=\beta_{i+1}$.
\medskip

According to the above, there is at least one real root of $H_1(f)$ in each of Int$(\beta_i,\beta_{i+1})$, $i \in \{1,\ldots,d-1\}$, where $\gamma_1\in\text{Int}(\beta_1,\beta_2)$, \ldots, $\gamma_{m_1-1}\in\text{Int}(\beta_{m_1-1},\beta_{m_1})$ are the first $m_1-1$  roots of $H_1(f)$ (in fact, the same root counted with multiplicity $m_1-1$) and similarly for the other multiple roots of $f$.

We have accounted for a total of $s-1+(m_1-1)+\cdots+(m_s-1)=m_1+\cdots+m_s-1 = d-1$ real roots of $H_1(f)$ counted with multiplicities. Hence $H_1(f)$ has no roots, real or complex, other than the ones listed above, and the result follows.
$\Box$

\begin{cor} Let $f$ be a polynomial of degree $d$ with $d$ real roots, counted with multipliciites. Then each of the $H_i(f)$,
$i\in\{1,\ldots,d-1\}$, has $d-i$ real roots, counted with multiplicity. In other words, all the roots of $H_i(f)$ are real.  
\end{cor}

Next, we recall a result from \cite{DdJ} on almost counterexamples to the Casas-Alvero conjecture.

\begin{deft} Fix an $i \in \{1,\ldots,d-1\}$. An {\bf almost counterexample to the Casas-Alvero conjecture of level} $i$ is a polynomial $f$ that has a common root with $H_j(f)$ for all $j\in\{1,\ldots,d-1\}\setminus\{i\}$ but is not a power of a linear polynomial.
\end{deft}

\noi{\bf Notation.} Given a polynomial $f$ of degree $d$ with $d$ real roots, for a pair $(k, m)$ of integers with $1\le k\le d-1$ and $1 \le m \le d-k$ we write $\alpha_{k,m}(f)$ for the $m$-th root of $H_k(f)$, where the roots of $H_k(f)$ are ordered (weakly) increasingly.
\medskip

We state the next theorem in a somewhat stronger form than in \cite{DdJ}: the extra information about the recycled roots
$\alpha_{k_j,m_j}(f)$ does not appear explicitly in the statement of the result in \cite{DdJ}, but is shown in the course of its proof.

\begin{thm}[J. Draisma--J. P. de Jong \cite{DdJ}, Theorem 5]\label{thm:ddj}

\noi Fix $d-2$ pairs of integers
\[
(k_j,m_j),\quad j \in \{1,\ldots,d-2\},
\]
with
\[
1\le k_1<k_2<\dots<k_{d-2}\le d-1
\]
and $1\le m_j\le d-k_j$. There exists a polynomial $f \in \R[x]$ with $f(0)=f(1)=0$, all of whose roots are real and lie in $[0,1]$, such that $\alpha_{k_j,m_j}(f)$ is a root of $f$ for all $j \in \{1,\ldots,d-2\}$ (in particular, $f$ is an almost counterexample to the Casas-Alvero conjecture of level $i$, where $i$ is the unique element of the set
$\{1,\dots,d-1\}\setminus\{k_1,\dots,k_{d-2}\}$). 
\end{thm}

We also recall the following result, Theorem 13 of \cite{CLO}:

\begin{thm} \label{thm:clo}
Assume that $f$ is a counterexample to the Casas-Alvero Conjecture. Then $f$ has at least five distinct roots.
\end{thm}

\noi Proof of Theorem \ref{thm:princ}: We argue by contradiction. Assume that
\begin{equation} \label{eq:rac}
R_i \in \sqrt{(R_1,\ldots,\breve{R_i},\ldots,R_{d-1})}.
\end{equation}
Let $f$ be an almost counterexample of level $i$ to the Casas-Alvero conjecture with $m_j=1$ for all $j \in \{1,\ldots,d-2\}$, given by Theorem \ref{thm:ddj}. By (\ref{eq:rac}), $f$ is a Casas-Alvero polynomial that is a counterexample to the Casas-Alvero conjecture. By  definition of $f$, $\alpha_{k_j1}(f)$ is the first root of $H_{k_j}(f)$ and also a root of $f$, for all $j \in \{1,\ldots,d-2\}$. In particular, since $i\in\{d-3,d-2,d-1\}$, $\alpha_{\ell 1}(f)$ is the first root of $H_\ell(f)$ and a root of $f$ for all $\ell\in\{1,\ldots,d-4\}$.      

By Proposition \ref{prop:sturm}, $\alpha_{11}(f)$ is also the first root of $f$ (so $\alpha_{11}(f)=0$). Let $m$ denote the multiplicity of this root of $f$. Since 0 is a root of $f$ of multiplicity $m$, it is also the first root of $H_\ell(f)$ for all
$1\leq \ell\leq m-1$ (again by Proposition \ref{prop:sturm}). In other words, $\alpha_{\ell1}(f)=0$ for all $1\leq\ell\leq m-1$. By Theorem \ref{thm:clo}, $f$ has at least 5 distinct roots, hence $m\le d-4$. Let $\beta$ denote the first strictly positive root of $f$.

By Proposition \ref{prop:sturm}, there is a unique root $\beta^{(1)}$ of $H_1(f)$ with $\beta^{(1)}\in]0,\beta[$. If $m>2$, again by Proposition \ref{prop:sturm}, there is a unique root $\beta^{(2)}$ of $H_2(f)$ with
$\beta^{(2)}\in]0,\beta^{(1)}[\subset ]0,\beta[$. We continue like this recursively until $H_{m-1}(f)$, to show that there is a unique root $\beta^{(m-1)}$ of $H_{m-1}(f)$ with $\beta^{(m-1)}\in]0,\beta^{(m-2)}[ \subset ]0,\beta[$. Now, $H_{m}(f)$ has no root at 0, hence its first root $\alpha_{m1}(f)$ belongs to the open interval $]0,\beta^{(m-1)}[\subset]0,\beta[$. 

This contradicts the fact that $\alpha_{m1}(f)$ is a root of $f$. $\Box$

\end{document}